\documentclass[12pt,amsmath]{article}
\usepackage{amsfonts}

\begin{document}

\newsymbol\rtimes 226F
\newsymbol\ltimes 226E
\newcommand{\text}[1]{\mbox{{\rm #1}}}
\newcommand{\Rep}{\text{Rep}}
\newcommand{\Corep}{\text{Corep}}
\newcommand{\End}{\text{End}}
\newcommand{\Hom}{\text{Hom}}
\newcommand{\Irr}{\text{Irr}}
\newcommand{\Ad}{\text{Ad}}
\newcommand{\tr}{\text{tr}}
\newcommand{\length}{\text{length}}
\newcommand{\ot}{\otimes}
\newcommand{\qed}{\kern 5pt\vrule height8pt width6.5pt depth2pt}
\def\newtheorems{\newtheorem{theorem}{Theorem}[subsection]
                 \newtheorem{cor}[theorem]{Corollary}
                 \newtheorem{prop}[theorem]{Proposition}
                 \newtheorem{lemma}[theorem]{Lemma}
                 \newtheorem{defn}[theorem]{Definition}
                 \newtheorem{Theorem}{Theorem}[section]
                 \newtheorem{Corollary}[Theorem]{Corollary}
                 \newtheorem{Proposition}[Theorem]{Proposition}
                 \newtheorem{Lemma}[Theorem]{Lemma}
                 \newtheorem{Definition}[Theorem]{Definition}
                 \newtheorem{Example}[Theorem]{Example}
                 \newtheorem{Remark}[Theorem]{Remark}
                 \newtheorem{claim}[theorem]{Claim}
                 \newtheorem{sublemma}[theorem]{Sublemma}
                 \newtheorem{example}[theorem]{Example}
                 \newtheorem{remark}[theorem]{Remark}
                 \newtheorem{question}[theorem]{Question}
                 \newtheorem{Question}[Theorem]{Question}
                 \newtheorem{conjecture}{Conjecture}[subsection]}
\newtheorems
\newcommand{\proof}{\par\noindent{\bf Proof:}\quad}

\title{{\bf The classification of finite-dimensional triangular
Hopf algebras over an algebraically closed field of characteristic
$0$}}
\author{Pavel Etingof\\
Massachusetts Institute of Technology\\ Department of
Mathematics\\ Cambridge, MA 02139, USA\\ {\rm email:
etingof@math.mit.edu} \and Shlomo Gelaki\\ Technion-Israel
Institute of Technology\\ Department of Mathematics\\ Haifa 32000,
Israel\\ {\rm email: gelaki@math.technion.ac.il} }
\maketitle

\section{Introduction}

The first step towards the classification of finite-dimensional
triangular Hopf algebras $H$ over an algebraically closed field
$k$ of characteristic $0$ was taken in [EG1, Theorem 2.1] where it
was proved that if $H$ is {\em semisimple} then it is obtained
from the group algebra $k[G]$ of a finite group $G$ by twisting
its comultiplication in the sense of Drinfeld [Dr]. The proof of
this theorem relies in an essential way on a theorem of Deligne on
Tannakian categories [De1] which characterizes symmetric tensor
rigid categories over $k$ which are equivalent to representation
categories of affine proalgebraic groups intrinsically as those
categories in which the categorical dimensions of objects are
non-negative integers. One can apply Deligne Theorem [De1] since
the representation category $\Rep(H)$ of $H$ has this property
(maybe after modifying its commutativity constraint). Later on we
used [EG1, Theorem 2.1] and the theory of Movshev on twisting in
finite groups [M] to completely classify semisimple triangular
Hopf algebras over $k$ in terms of certain quadruples $(G,A,V,u)$
of group-theoretical data [EG2], and obtained a similar
classification in positive characteristic [EG2].

However, for {\em non-semisimple} finite-dimensional triangular
Hopf algebras $H$ over $k$ it is no longer true that the
categorical dimensions of objects in $\Rep(H)$ are non-negative
integers; so Deligne Theorem [De1] cannot be applied.
Nevertheless, in [AEG] it was realized that all known examples of
finite-dimensional triangular Hopf algebras $H$ over $k$ have the
{\em Chevalley property}; namely, the semisimple part of $H$ is
itself a Hopf algebra. Then in [AEG, Theorem 5.1.1], using [EG1,
Theorem 2.1] and hence Deligne Theorem [De1], it was proved that
$H$ has the Chevalley property if and only if it is obtained from
a certain modification of the supergroup algebra $k[G]$ of a
finite supergroup $G$ by twisting its comultiplication. Later on
we used [AEG, Theorem 5.1.1] to completely classify
finite-dimensional triangular Hopf algebras over $k$ with the
Chevalley property in terms of certain septuples $(G,W,A,Y,B,V,u)$
of group-theoretical data [EG3]. Nevertheless, we could not show
that this list contains all possible finite-dimensional triangular
Hopf algebras over $k$.

Very recently, Deligne has generalized his theorem on Tannakian
categories [De1] to the super-case [De2]. The purpose of this note
is to explain how this remarkable generalized theorem, combined
with the results from [AEG,EG1,EG2,EG3], lead to the complete and
explicit classification of finite-dimensional triangular Hopf
algebras over an algebraically closed field $k$ of characteristic
$0$, and to answer some questions from [AEG,G2] about triangular
Hopf algebras and symmetric rigid tensor categories over $k$,
positively. We also use Deligne Theorem and the lifting theorem of
Etingof-Nikshych-Ostrik [ENO] to prove that a symmetric fusion
category over a field of characteristic $p>0$, whose global
dimension is nonzero, is equivalent to a representation category
of a unique finite group whose order is not divisible by $p$ (see
Theorem \ref{pchar} below). We emphasize that although the results
of this note are useful for Hopf algebra theory, they (with the
exception of Theorem \ref{pchar}) follow from [De2] and
[AEG,EG1,EG2,EG3] in a fairly straightforward manner.

\noindent {\bf Acknowledgments} We are grateful to
A. Braverman, P. Deligne, D. Kazhdan, and V. Ostrik
for useful discussions.

\section{Preliminaries}
Throughout this note, unless otherwise stated, $k$ will denote an
algebraically closed field of characteristic $0$.

Recall that an {\em affine algebraic supergroup} is the spectrum
of a finitely generated supercommutative Hopf superalgebra over
$k$ (see [De2]). In other words, it is a functor $G$ from the
category of supercommutative algebras to the category of groups
defined by $A\mapsto G(A):=\Hom(H,A)$, where $\Hom(H,A)$ is the
set of algebra maps $H\to A$ preserving the parity (so called
functor of points). An inverse limit of affine algebraic
supergroups is called an {\em affine proalgebraic supergroup}. A
{\em finite supergroup} is an affine algebraic supergroup, whose
function algebra is finite-dimensional. In this case, the dual of
this function algebra (called the {\em supergroup algebra}) is a
supercocommutative Hopf superalgebra of the form $k[G\ltimes
V]=k[G]\ltimes \Lambda V$, where $G$ is a finite group and $V$ is
a finite-dimensional $k-$representation of $G$ (see e.g. [AEG] for
more details).

Let $G$ be an affine proalgebraic supergroup over $k$ and let $p:
\mathbb{Z}_2\to G$ be a morphism such that $\Ad(p(-1))$ is the
parity automorphism of $G$. Let $\Rep(G,p)$ denote the category of
all finite-dimensional algebraic representations of $G$ over $k$
in which $p(-1)$ acts as the parity operator. Then $\Rep(G,p)$ is
a $k-$linear abelian symmetric rigid tensor category with
$\End({\bf 1})=k,$ where ${\bf 1}$ denotes the unit object of
$\Rep(G,p)$ (see [DM]).

\begin{Definition} Let $\mathcal{C}$ be a $k-$linear (abelian)
symmetric tensor category which is equivalent to
$\Rep(G,p)$ for some $G,p$. Then $\mathcal{C}$ is said to be of
supergroup type. If in addition $G$ is a finite supergroup then
$\mathcal{C}$ is said to be of finite supergroup type.
\end{Definition}

\section{Deligne Theorem}

Let $\mathcal{C}$ be a $k-$linear abelian symmetric rigid tensor
category with $\End({\bf 1})=k$. Recall that for any $X$ in
$\mathcal{C}$ its length, denoted by $\length(X)$, is defined to
be the maximal possible length of a strictly increasing filtration
of $X$. We will always assume that all objects in $\mathcal{C}$
have finite length. We are now ready to state Deligne Theorem.

\begin{Theorem}\label{c2} {\rm (see [De2],
Proposition 0.5, Theorem 0.6, and  the sentence after Theorem
0.6)} Suppose that for any object $X$ in $\mathcal{C}$ there
exists a constant $d(X)>0$ such that $\length(X^{\otimes n})$ is
dominated by $d(X)^n$. Then $\mathcal{C}$ is of supergroup type.
\end{Theorem}

\begin{Corollary}\label{c3} {\rm ([De2], Corollaries 0.7 and 0.8)}
If $\mathcal{C}$ has finitely many classes of simple objects,
 then $\mathcal{C}$ is of
supergroup type. In particular, if in addition $\mathcal{C}$ is
semisimple then $\mathcal{C}$ is equivalent to $\Rep(G)$ where $G$
is a finite group, possibly with a modified symmetric structure.
\end{Corollary}

If $\mathcal{C}$ is equivalent, as a $k-$linear abelian category,
to $\Rep(A)$, where $A$ is a finite-dimensional algebra, then
$\mathcal{C}$ is said to be {\em finite}. It is known that this
condition is equivalent to the condition that $\mathcal{C}$ has
finitely many isomorphism classes of simple objects, and any
simple object has a projective cover. It follows from Corollary
\ref{c3} that if $\mathcal C$ is finite then it is of finite
supergroup type.

\section{Applications to Hopf algebras}

Corollary \ref{c3} implies Theorem 2.1 of [EG1] on the
classification of triangular semisimple Hopf algebras.

\begin{Theorem}\label{main1} {\rm [EG1, Theorem 2.1]}
Let $(H,R)$ be a semisimple triangular Hopf algebra over $k$, with
Drinfeld element $u$. Set $R_u:=1/2(1\ot 1+1\ot u+u\ot 1-u\ot u)$
and $\tilde R:=RR_u$. Then there exist a finite group $G$ and a
twist $J\in k[G]\ot k[G]$ such that $(H,\tilde R)$ and
$(k[G]^J,J_{21}^{-1}J)$ are isomorphic as triangular Hopf
algebras.
\end{Theorem}

\proof Let $\Rep(H)$ denote the $k-$linear abelian symmetric
tensor rigid category of all finite-dimensional
$k-$representations of $H$. Clearly, $\End({\bf 1})=k$. By
Corollary \ref{c3}, there exists a finite group $G$ such that
$\Rep(H)$ is equivalent to $\Rep(G)$, possibly with modified
symmetric structure. The rest is as in [EG1] (see also [G1]). \qed

\begin{Remark} {\rm The proof we gave in [EG1] relied on the weaker
version of Deligne Theorem [De1] and hence required some Hopf
algebra theory e.g. Larson-Radford Theorem that the antipode of a
semisimple Hopf algebra over $k$ is an involution [LR]. We stress
that Hopf algebra theory is no longer needed for the proof of
Theorem \ref{main1}.}
\end{Remark}

Recall [AEG] that a finite-dimensional triangular Hopf algebra
$(H,R)$ is called a {\em modified supergroup algebra} if its
$R-$matrix $R$ is of rank $\le 2$. The reason for this terminology
can be found in Corollaries 2.3.5 and 3.3.3 in [AEG] where it is
proved that such finite-dimensional triangular Hopf algebras correspond
to (finite) supergroup algebras. The correspondence respects the
tensor categories of representations [AEG, Theorem 3.1.1] and the
twisting procedure [AEG, Proposition 3.2.1]. Recall also that $H$
is said to have the {\em Chevalley property} if its quotient by
its radical is a Hopf algebra itself [AEG]. In [AEG, Theorem
5.1.1] it is proved that $H$ is twist equivalent to a modified
supergroup algebra (by twisting of comultiplication) if and only
if $H$ has the Chevalley property. In Question 5.5.1 in [AEG] we
asked if any finite-dimensional triangular Hopf algebra over $k$
has the Chevalley property. We now have

\begin{Theorem}\label{main2}
Let $H$ be a finite-dimensional triangular Hopf algebra over $k$.
Then $H$ is twist equivalent to a modified supergroup algebra. In
particular, $H$ has the Chevalley property.
\end{Theorem}

\proof This follows from Corollary \ref{c3} and the preceding
remarks. \qed

Recall from [EG3] that a {\em triangular septuple} is a septuple
$(G,W,A,Y,B,V,u)$ where $G$ is a finite group, $W$ is a
finite-dimensional $k-$representation of $G$, $A$ is a subgroup of
$G$, $Y$ is an $A-$invariant subspace of $W$, $B$ is an
$A-$invariant nondegenerate element in $S^2Y$, $V$ is an
irreducible projective $k-$representation of $A$ of dimension
$|A|^{1/2}$, and $u\in G$ is a central element of order $\le 2$
acting by $-1$ on $W$. In [EG3, Section 2] it is explained how to
assign a finite-dimensional triangular Hopf algebra with the
Chevalley property $H(G,W,A,Y,B,V,u)$ over $k$ to any triangular
septuple. Therefore, Theorem \ref{main2} and [EG3, Theorem 2.2]
imply now the following explicit classification of
finite-dimensional triangular Hopf algebras over $k$.

\begin{Theorem}\label{main3}
The assignment $(G,W,A,Y,B,V,u)\to H(G,W,A,Y,B,V,u)$ is a
bijection between:
\begin{enumerate}

\item isomorphism classes of triangular septuples, and

\item isomorphism classes of finite-dimensional triangular Hopf
algebras over $k$.
\end{enumerate}
\end{Theorem}

\begin{Remark} {\rm In [EG2, Theorem 5.2] we proved that for
$W=Y=B=0$ the bijection given in Theorem \ref{main3} reduces to a
bijection between quadruples $(G,A,V,u)$ and semisimple triangular
Hopf algebras over $k$.}
\end{Remark}

As a result of Theorem \ref{main2} and [AEG, Remark 5.5.2] we can
now answer Question 3.6 from [G2] positively.

\begin{Theorem}\label{main4}
Let $H$ be a finite-dimensional triangular Hopf algebra over $k$
and let $u$ be its Drinfeld element. Then $u^2=1$ and consequently
$S^4=id$. Moreover, if $\dim(H)$ is odd then $u=1$ and $H$ is
semisimple.
\end{Theorem}

Finally, let us explain how Deligne Theorem can be applied to
cotriangular Hopf algebras $H$ over $k$. Let $\Corep(H)$ denote
the $k-$linear abelian symmetric tensor rigid category of all
finite-dimensional $k-$corepresentations of $H$. Clearly,
$\End({\bf 1})=k$.

\begin{Theorem}\label{main5}
Let $H$ be a cotriangular Hopf algebra over $k$. Then the category
$\Corep(H)$ is equivalent to the category $\Rep(G,p)$ for a unique
affine proalgebraic supergroup $G$ and morphism $p$.
\end{Theorem}

\proof Set $\mathcal{C}:=\Corep(H)$. If $X$ is in $\mathcal{C}$
then it is clear that $\length(X^{\otimes n})$ is at most
$\dim(X)^n$, where $\dim(X)$ is the usual linear algebraic
dimension of $X$. Thus, by Theorem \ref{c2}, $\mathcal{C}$ is of
the form $\Rep(G,p)$. The uniqueness of $(G,p)$ follows from the
uniqueness of the super fiber functor, which follows from Sections
3 and 4 of [De2]. \qed

\bigskip
\noindent {\bf Important remark.} In [EG4, Theorem 3.3] we showed
that any {\em pseudoinvolutive} cotriangular Hopf algebra over $k$
(i.e. such that $\tr(S^2_{|C})=\dim(C)$ on any finite-dimensional
subcoalgebra $C$) is twist equivalent (by twisting of
multiplication) to $\mathcal{O}(G)$ for an affine proalgebraic
{\em group} $G$, and vice versa. One might expect that this
correspondence would extend to supergroups, if one drops the
pseudoinvolutivity condition. As we saw, this is definitely true
in the finite-dimensional case. Nevertheless, in the
infinite-dimensional case, such a generalization fails, and the
situation is much more nontrivial. Namely, Theorem \ref{main5}
implies that the coalgebras $H$ and $\mathcal{O}(G)$ are Morita
equivalent, but it does {\bf not} imply that they are isomorphic,
since the equivalence of Theorem \ref{main5} need not preserve
linear algebraic dimensions (as, unlike in [EG4], they need not be
equal to the categorical dimensions). In fact, even for $G=SL(2)$,
for any integer $N>2$ there exist cotriangular Hopf algebras $H$
with ${\rm Comod}(H)=\Rep(G)$, such that the 2-dimensional vector
representation of $G$ corresponds to an N-dimensional object in
${\rm Comod}(H)$. (For examples of such Hopf algebras, see
[GM,B]). This shows that the theory of (infinite-dimensional)
cotriangular Hopf algebras is much richer than the theory of
triangular Hopf algebras.

\section{Applications to tensor categories}

We first use Deligne Theorem to answer Question 5.5.5 in [AEG].
Recall from [AEG] that a tensor category $\mathcal{C}$ is said to
have the Chevalley property if the tensor product of any two
simple objects in $\mathcal{C}$ is semisimple.

\begin{Theorem} \label{main6} Let $\mathcal{C}$ be a finite symmetric
tensor rigid category with $\End({\bf 1})=k$. Then $\mathcal{C}$
has the Chevalley property.
\end{Theorem}

\proof This follows from Corollary \ref{c3} since $\Rep(G,p)$ has
the Chevalley property when $G$ is a finite supergroup. \qed

\begin{Theorem} \label{main7} Let $\mathcal{C}$ be a finite symmetric
tensor rigid category with $\End({\bf 1})=k$. Then $\mathcal{C}$
is equivalent to a category $\Rep(H)$ where $H$ is a
finite-dimensional triangular Hopf algebra with $R-$matrix of rank
$\le 2$.
\end{Theorem}

\proof By Corollary \ref{c3}, $\mathcal{C}$ is equivalent to
the category of representations of a finite supergroup $G\ltimes
V$ on supervector spaces, in which a fixed element of order $2$
acting by parity on $G\ltimes V$ is represented by the parity
operator. Thus, $\mathcal{C}$ is equivalent to the category of
representations of the cocommutative triangular Hopf superalgebra
$k[G\ltimes V]$. Modifying it into a finite-dimensional triangular
Hopf algebra $H$ with $R-$matrix of rank $\le 2$ (as in [AEG,
Corollary 3.3.3]), we obtain the theorem. \qed

\begin{Remark} {\rm In Questions 5.5.5 and 5.5.6 in [AEG] we accidentally
omitted the assumption that the category $\mathcal{C}$ has enough
projectives (i.e. each simple object has a projective cover).
Without this condition one may take for example the
category $\mathcal{C}:=\Rep(G_a)$ of algebraic representations of
the additive group; it is definitely not a representation category of a
finite-dimensional triangular Hopf algebra.}
\end{Remark}

We now apply Deligne Theorem to symmetric fusion categories over
fields with positive characteristics. 
Let $k$ be an algebraically closed field of characteristic $p>0$.
Recall that if $\mathcal{C}$
is a fusion category ($=$ semisimple rigid tensor category with
finitely many simple objects and $\End({\bf 1})=k$), the global
dimension $\dim(\mathcal{C})$ of $\mathcal{C}$ is defined to be
the sum of squares of dimensions of its simple objects (see e.g.
[ENO]).

\begin{Theorem}\label{pchar} Let $\mathcal{C}$ be a symmetric
fusion category over a field $k$ of positive characteristic $p$.
Assume that $\dim(\mathcal{C})$ is nonzero. Then $\mathcal{C}$ is
equivalent to $\Rep_k(G)$ for a unique (up to isomorphism) finite
group $G$ whose order is not divisible by $p$.
\end{Theorem}

\proof According to [ENO, Theorem 9.3], the category $\mathcal{C}$
admits a unique lifting to a category $\mathcal{C}_{W(k)}$ over
the ring of Witt vectors $W(k)$. Let $K$ be the algebraic closure
of the field of quotients $Q$ of $W(k)$. Let $\mathcal{C}_Q$,
$\mathcal{C}_K$ be the categories obtained from
$\mathcal{C}_{W(k)}$ by extension of scalars to $Q,K$
respectively. By Corollary \ref{c3}, there exists a finite group
$G$ such that $\mathcal{C}_K=\Rep_K(G)$. Thus there exists a
finite extension $Q'$ of $Q$ such that
$\mathcal{C}_{Q'}=\Rep_{Q'}(G)$.

The order of the group $G$ is the global dimension of
$\mathcal{C}_K$, hence it is equal to the global dimension of
$\mathcal{C}$ modulo $p$. Thus, $p$ does not divide $|G|$.

Let $W'$ be the ring of integers in $Q'$. It is a local ring. Let
$I$ be the maximal ideal in $W'$. Then the residue field $W'/I$ is
equal to $k$.

Consider the tensor categories $\mathcal{C}_{W'}$ and
$\Rep_{W'}(G)$. The reductions modulo $I$ of these categories
(namely, $\mathcal{C}$ and $\Rep_k(G)$, respectively) have nonzero
global dimensions. The localizations $\mathcal{C}_{Q'}$,
$\Rep_{W'}(G)$ are equivalent. By [ENO, Theorem 9.6(ii)], this implies
that the reductions $\mathcal{C}$ and $\Rep_k(G)$ are equivalent
as well, as symmetric categories.

The uniqueness of $G$ is well known (see [DM]). The theorem is
proved. \qed

\begin{Remark} {\rm if $\dim(\mathcal{C})=0$, the theorem is false, and
much more interesting categories than $\Rep_k(G)$ can occur. A
counterexample is the reduction mod $p$ of the fusion category of
representations of $U_q(sl_2)$, $q=e^{\pi i/p}$. This category is
symmetric, but is not equivalent to the category of
representations of a finite group, since the Frobenius-Perron
dimensions of its objects are not integers.}
(see \cite{eno}, Remark 9.5). 
\end{Remark}

\end{document}